\documentclass[12pt]{amsart}
\usepackage[a4paper]{geometry}               
\usepackage{graphicx}
\usepackage{amssymb}
\usepackage{amsthm}
\usepackage{filecontents}

\def\bstates{\varOmega}
\def\ISstates{\varSigma}
\def\Ex{\mathop{\null\mathbb{E}}\nolimits}
\def\BIS{\mathsf{\#BIS}}
\def\numP{\mathsf{\#P}}
\def\NP{\mathsf{NP}}

\def\bpi{\pi_\bstates}
\def\omega{\pi_\ISstates}

\newtheorem{theorem}{Theorem}
\newtheorem{lemma}[theorem]{Lemma}

\makeatletter
\def\prob#1#2#3{\goodbreak\begin{list}{}{\labelwidth\z@ \itemindent-\leftmargin
                        \itemsep\z@  \topsep6\p@\@plus6\p@
                        \let\makelabel\descriptionlabel}
                \item[\it Name]#1
                \item[\it Instance]#2
                \item[\it Output]#3
                \end{list}}
\makeatother

\title[The Ge-\v Stefankovi\v c Process]{A Counterexample to rapid mixing\\ 
of the Ge-\v Stefankovi\v c Process}
\author{Leslie Ann Goldberg}
\author{Mark Jerrum}

\address{Leslie Ann Goldberg\\
Department of Computer Science\\
University of Liverpool\\Ashton Building\\
Liverpool L69 3BX\\
United Kingdom.}

\address{Mark Jerrum\\
School of Mathematical Sciences\\
Queen Mary, University of London\\
Mile End Road\\
London E1 NS\\
United Kingdom.}

\thanks{The work described in this paper was partly supported by EPSRC Research Grant
(refs\ EP/I011528/1,  EP/I011935/1 and EP/I012087/1) ``Computational Counting''.}

\begin{document}
\maketitle

\begin{abstract}
Ge and \v Stefankovi\v c have recently introduced a novel two-variable graph polynomial.
When specialised to a bipartite graphs~$G$ and evaluated at the point $(\frac12,1)$ 
this polynomial gives 
the number of independent sets in the graph.  Inspired by this polynomial, they also
introduced a Markov chain which, if rapidly mixing, would provide an efficient
sampling procedure for independent sets in~$G$.  This sampling procedure 
in turn would imply the existence of efficient approximation algorithms for a number of 
significant counting problems whose complexity is so far unresolved.
The proposed Markov chain is promising, in the sense that it overcomes the most obvious 
barrier to mixing.  However, we show here, by exhibiting a sequence of counterexamples, 
that the mixing time of their Markov chain 
is exponential in the size of the input when the input is chosen from a particular infinite family of bipartite graphs. 
\end{abstract}

\section{Overview}\label{sec:overview}
Consider the following basic computational problem:
\prob{$\BIS$.}{A bipartite graph $G$.}{The number of independent sets
in~$G$.}
It has long been know that $\BIS$ is $\numP$-complete~\cite{PB}, and hence presumably
intractable, if we insist on an exact solution.  However, the computational complexity of 
approximating~$\BIS$ remains a fascinating open problem.  The standard notion of 
efficient approximation in the context of counting problems is the ``fully polynomial
approximation scheme'' or FPRAS\null.  Roughly speaking, an FPRAS is a polynomial-time
randomised algorithm that produces an estimate that is close in relative error to the true solution 
with high probability.  (See \cite[Defn~11.2]{MotRag} for a precise definition.)
The most satisfactory situation would be either to have an FPRAS for $\BIS$, or a proof that
$\BIS$ is $\NP$-hard to approximate.  However, neither of situations is known to 
occur.

Dyer
et al.~\cite{DGGJ} noted that a number of counting problems are equivalent 
to $\BIS$ under approximation-preserving reducibility, and further
$\BIS$-equivalent problems have been presented in subsequent work~\cite{StableMatchings,GJ07}.  
Since no FPRAS has been found for any of the counting 
problems in this equivalence class, it is becoming standard to progress under the assumption
that $\BIS$ (and hence 
each of
the $\BIS$ equivalent problems) does not admit an FPRAS\null.
So finding an FPRAS for $\BIS$ at this stage would be a significant development.  Not
only would it imply the existence of an FPRAS for several natural counting problems ---
such as counting downsets in a partial order, or evaluating the partition function of
the ferromagnetic Ising model with local fields --- but it would also resolve the 
complexity of approximating $\BIS$ in the opposite direction to the one many people expect.

The most fruitful approach to designing efficient approximation algorithms for counting problems
has been Markov chain Monte Carlo (MCMC)\null.  A direct application of MCMC to $\BIS$ would work as follows.
Given a bipartite graph~$G$ with $n$~vertices, 
consider the Markov chain whose 
state space, $\Omega$, is the set of all independent sets in~$G$,
and whose transition probabilities $P(\cdot,\cdot)$ are as follows,
where $\oplus$ denotes symmetric difference and $H(I) = \{I' \in \Omega \mid |I \oplus I'| = 1\}$:
$$
P(I,I')= \begin{cases}
1-H(I)/2n, & \text{if $I'=I$};\\
1/2n, &\text{if  $I'\in H(I)$};\\
0, &\text{otherwise}.  
\end{cases}
$$ 
It is easy to check that this Markov chain has the uniform distribution on independent
sets as its unique stationary distribution.  So, simulating the Markov 
chain for sufficiently many steps would enable us to sample independent sets
nearly uniformly.  From there it is a short step to estimating
the number of independent sets~\cite[\S3.2]{ETH}.  

To obtain an FPRAS from this approach, one requires that the 
Markov chain on independent sets is rapidly mixing, i.e., that it is close to
the stationary distribution after a number of steps that is polynomial in~$n$.
Unfortunately, it is clear that the proposed Markov chain does not have this
property.  Consider the complete bipartite graph with equal numbers of vertices 
in the left and right blocks of the bipartition.  There are $2^{n/2}-1$ 
independent sets that have non-empty intersection with the left block, 
and the same number with non-empty intersection with the right.  Any sequence
of transitions which starts in a left-oriented independent set and ends in a 
right-oriented one must necessarily pass through the empty independent set.
Informally, the empty independent set presents an obstruction to rapid mixing
by forming a constriction in the state space.  This intuition can be made 
rigorous by noting that the ``conductance'' of the Markov chain is exponentially
small, which implies exponential (in~$n$) mixing time \cite[Claim~2.3]{DFJ}.
In fact, it is not even necessary to have 
a dense graph in order to obtain
such a constriction:  degree~6 will do \cite[Thm 2.1]{DFJ}.

Ge and \v Stefankovi\v c~\cite{GSconf} have recently introduced an intriguing graph 
polynomial $R_2'(G;\lambda,\mu)$, in two indeterminates $\lambda$ and~$\mu$,
that is associated with a bipartite graph~$G$.  At the point $(\lambda,\mu)=(\frac12,1)$
it counts independent set in~$G$;  specifically,
the number of independent sets in~$G$ is given by $2^{n-m}R_2'(G;\frac12,1)$,
where $n$ is the number of vertices, and $m$ the number of edges in~$G$  
\cite[Thm~4]{GSconf}.
This polynomial inspires them to propose a new Markov chain 
\cite[Defn~6]{GSconf}
that potentially could be used to sample independent sets from 
a bipartite graph and hence provide an approximation algorithms for~$\BIS$.  
The Markov chain, which is described below,  
is very different from the one discussed earlier.  In particular, its states
are subsets of the edge set of~$G$ rather than subsets of the vertex set.
Thus, sampling an independent set of~$G$ is a two-stage procedure:  (a)~sample
an edge subset~$A$ of~$G$ from the appropriate distribution, and then (b)~sample
an independent set from a distribution conditioned on~$A$.  Details will be given 
below.

The encouraging aspect of this new Markov chain, which we call the
Ge-\v Stefankovi\v c Process, or GS Process for short, 
is that it is immune to the obvious counterexamples,
such as the complete bipartite graph.  
Unfortunately, with a certain amount of effort it is possible to
find a counterexample to rapid mixing.  In the following section we describe 
the GS Process and construct a sequence of graphs on which its mixing time is
exponential (in the number of vertices of the graph).  Although this counterexample 
rules out their Markov chain as an approach to constructing a general FPRAS for~$\BIS$,  
we may still hope that it provides an efficient algorithm for some restricted class
of graphs.
For example, \cite[Theorem 7]{GSconf} shows that it provides an efficient algorithm on trees.

\section{The Ge-\v Stefankovi\v c Process}

Before stating our result, we need to formalise what we mean by
mixing, rapid or otherwise.  Let $(X_t)$ be an
ergodic Markov chain with state
space $\bstates$, distribution~$p_t$ at time $t$, 
and unique stationary distribution~$\pi$.  
Let $x_0\in\bstates$ be the initial state of the chain,
so that $p_0$ assigns unit mass to state~$x_0$.
Define the {\it mixing time $\tau(x_0)$
with initial state~$x_0\in\bstates$}, 
as the first time~$t$ at which
$\frac12{\|p_t-\pi\|}_1 \leq e^{-1}$, i.e., at which the distance between the $t$-step
and stationary distributions as at most $e^{-1}$ in total variation;
then define the {\it mixing time}~$\tau$ as the maximum of 
$\tau(x_0)$ over all choices of initial state~$x_0$.

Suppose $G=(U\cup V,E)$ is a bipartite graph, where $U,V$ are disjoint sets forming the
vertex bipartition, and $E$ is the edge set.  We are interested in two probability
spaces, $(\bstates,\bpi)$ and $(\ISstates,\omega)$, where $\bstates=2^E$ and $\ISstates=2^U$.
We construct the probability distributions $\bpi:\bstates\to[0,1]$ and $\omega:\ISstates\to[0,1]$
with the help of a certain
consistency relation~$\chi$ on $\ISstates\times\bstates$, which is defined as follows.
For a pair $(I,A)\in\ISstates\times\bstates$, consider the 
subgraph of $(U\cup V,A)$ induced by the vertex set $I\cup V$.  
We say that the relation $\chi(I,A)$ holds iff every vertex of~$V$ 
has even degree in this subgraph.  
Start with the probability space of consistent pairs $\{(I,A)\in\ISstates\times\bstates\mid
\chi(I,A)\}$ with the uniform distribution.  
Then $\bpi$ (respectively~$\omega$) is the induced marginal distribution on~$\bstates$
(respectively~$\ISstates$).   We call $\omega$ the {\it marginal BIS distribution\/}
on~$\ISstates$.  It is shown in 
\cite[Lemma~10]{GSconf}
that $\omega$ is also the distribution induced on~$U$
by a uniform random independent set in~$G$, justifying the name.
In~\cite{GSconf}, $\bpi(A)$, for $A\in\bstates$ is defined in terms 
of the rank of~$A$, viewed as a bipartite adjacency matrix over $\mathbb{F}_2$;
this definition is equivalent to the one given here.

The GS-process is an ergodic ``single bond flip''
Markov chain on state space~$\bstates$ which has 
stationary distribution~$\bpi$.   The exact definition of this Markov chain 
is not important to us, as our counterexample applies to any Markov chain
on state space $\bstates$ with stationary distribution~$\bpi$ 
that does not change too many edges in one step.  In order to formalise this
last requirement, say that a Markov chain 
with transition probabilities $P:\bstates^2\to[0,1]$,
is {\it $d$-cautious\/} if
$$
 P(A,A')>0\implies |A\oplus A'|\leq d,\quad\text{for all $A,A'\in\bstates$}.
$$
The GS Process is a $1$-cautious Markov chain.  Our negative result applies
to all $d$-cautious chains, where $d$ depends 
at most linearly on the number of vertices of~$G$.

\section{A counterexample to rapid mixing}

The following lemma 
(taken from \cite[Claim 2.3]{DFJ})  
packages the conductance argument 
in a convenient way for us to obtain explicit lower bounds
on mixing time.

\begin{lemma}\label{lem:torpid}
Let $(X_t)$ be a Markov chain with state
space $\bstates$, transition matrix~$P$
and stationary distribution~$\pi$. 
Let $\{S,T\}$ be a partition of $\bstates$ such that $\pi(S) \leq \frac12$,
and $C\subset\bstates$ be a set of states that form a 
``barrier'' in the sense that $P(s,t)=0$ whenever 
$s\in S\setminus C$ and $t\in T\setminus C$. 
Then the mixing time of the Markov chain is at least $\pi(S)/8\pi(C)$.
\end{lemma}

Let $n,m$ be positive integers such that $(3/2)^{m}\leq2^n-1<(3/2)^{m+1}$.  
Note that for every $n$ there is a unique~$m$ satisfying the inequalities, 
and that $m$ depends linearly on $n$, asymptotically.  The 
counterexample graph (actually sequence of graphs indexed by~$n$)
$G_n=(U'\cup V\cup U'',E)$ has vertex set $U'\cup V\cup U''$ where
$|U'|=n$ and $|V|=|U''|=m$.  The edge set is $E=U'\times V\cup M$, where $M$ is 
a perfect matching of the vertices in $V$ and $U''$. Thus, (a)~$G_n$ has bipartition
$(U,V)$ where $U=U'\cup U''$,  (b)~$U'$, $V$ and $U''$ are all independent sets,
(c)~the edges between $U'$ and $V$ form a complete graph, and (d)~the 
edges between $V$ and $U''$ form a matching.

Partition $\ISstates$ as $\ISstates=\ISstates_0\cup\ISstates_1$, where 
$\ISstates_0=\{I\in\ISstates\mid I\cap U'=\emptyset\}$ and $\ISstates_1=\ISstates\setminus
\ISstates_0$.  Observe there are $3^m$ independent sets in~$G_n$ that exclude all vertices
in~$U'$, and $(2^n-1)2^m$ that include some vertex.  Since $\omega$ is the 
marginal distribution of independent sets in~$G_n$,
$$
\omega(\ISstates_0)=\frac{3^m}{3^m+(2^n-1)2^m}=\alpha,
$$
where $\frac25<\alpha\leq\frac12$, by choice of $n,m$.
So $\ISstates_0\cup\ISstates_1$ is a nearly
balanced partition of the state-space $\ISstates$.
Also it is easy to check that the cut defined by this partition
is a witness to the conductance of the ``single site flip'' Markov chain
of \S\ref{sec:overview}
being exponentially small in~$n$.
This implies that the mixing time of
the single site flip Markov chain is exponential in~$n$ (which, of course, 
was never in doubt).  Next we identify a partition $\bstates_0\cup\bstates_1$
that mirrors the partition $\ISstates_0\cup\ISstates_1$, and itself witnesses
exponentially small conductance of the GS Process.

Define the {\it weight} $w(A)$ of $A\in\bstates$ to be $w(A)=|A\cap M|$.
Partition $\bstates$ as $\bstates=\bstates_0\cup\bstates_1$, where 
$\bstates_0=\{A\in\bstates\mid w(A)\leq \frac5{12}m\}$ 
and $\bstates_1=\bstates\setminus\bstates_0$.
We aim to show that the weights of states in $\bstates$ are concentrated around
$\frac13m$ and $\frac12m$, and there are exponentially few states near the 
boundary of $\bstates_0$ and $\bstates_1$.
With a view to applying Lemma~\ref{lem:torpid}, define a ``barrier set'' (of states)
by  
$$
C=\{A\in\bstates\mid \tfrac9{24}m\leq w(A)\leq \tfrac{11}{24}m\}.
$$
It is not clear how to sample a state~$A$ from the distribution $(\bstates,\bpi)$
directly, so instead we sample a state~$I$ from $(\ISstates,\omega)$ and then sample
u.a.r.\ a state $A$ consistent with $I$, i.e., satisfying $\chi(I,A)$.  This amounts
to the same thing.

Suppose we start with a state $I$ sampled from $(\ISstates,\omega)$, conditional
on $I\in\ISstates_0$.  The set $I\cap U''$ is determined by 
a Bernoulli process with success probability~$\frac13$.  (For each edge~$e$ in~$M$ 
there are three possibilities for the restriction of the independent set~$I$ to~$e$,
and only one of them includes a vertex from~$U''$.
These choices are independent for each $e\in M$.)  When we come to select a 
random consistent edge set~$A$, we must exclude all edges in~$M$ that are incident to 
a vertex in~$I\cap U''$.  The other edges in $M$ are free to be included 
or excluded.  So the set of edges $A\cap M$ is determined by a Bernoulli 
process with success probability $\frac13$.  Thus $\Ex(w(A))=\frac13m$ and, by 
a Chernoff bound,
$\Pr(w(A)\geq \frac{9}{24}m)$ is exponentially small in~$m$.
Specifically, 
\begin{equation}\label{eq:bd1}
\Pr(A\in C)\leq\Pr\big(w(S)\geq \tfrac{9}{24}m\big)\leq\exp(-m/576)
\end{equation}
by \cite[Thm~4.4(2)]{MitzUpfal}, setting $\delta=\frac18$ and $\mu=\frac13m$.

Now suppose $I$ is sampled from $(\ISstates,\omega)$, conditional
on $I\in\ISstates_1$.  
Now select a uniform random~$A$, conditional on the event $\chi(I,A)$.
We argue that the probability that a given edge~$e=(v,u)$ 
of~$M$  
is included in~$A$ is~$\frac12$, independent of all the other 
edges of~$M$.  
Suppose $v\in V$ and $(v,u)\in M$.
Imagine we are deciding which edges incident to~$v$ 
are to be included in~$A$.  First we decide whether to include the edge $(v,u)$
itself.  In selecting the remaining edges for~$A$ from the $n$ available, 
we just have to make sure that the parity of $A\cap(\{v\}\times (I\cap U'))$ is odd,
if $(v,u)\in A$ and $u\in I$, and even otherwise.  Since $I\cap U'\not=\emptyset$,
the number of ways to do this is $2^{n-1}$, independent of whether we included
edge~$e$ in the first place.  It follows that the
set of edges $A\cap M$ is determined by a Bernoulli 
process with success probability $\frac12$.  Thus $\Ex(w(A))=\frac12m$ and, by 
Chernoff, $\Pr(w(A)\leq \tfrac{11}{24}m)$ is exponentially small in~$n$.
Specifically, 
\begin{equation}\label{eq:bd2}
\Pr(A\in C)\leq\Pr\big(w(S)\leq \tfrac{11}{24}m\big)\leq\exp(-m/576) 
\end{equation}
by~\cite[Thm~4.5(2)]{MitzUpfal}, setting $\delta=\frac1{12}$ and $\mu=\frac12m$.

We see now that the partition $\bstates_0\cup\bstates_1=\bstates$ is balanced, since
$
\bpi(\bstates_0)
=\alpha\pm o(1)$ and $\frac25<\alpha\leq\frac12$.
Moreover, from (\ref{eq:bd1}) and~(\ref{eq:bd2}), 
$\Pr(\frac9{24}\leq w(A)\leq \frac{11}{24})$ is exponentially small
when $A$ is selected from the distribution $(\bstates,\bpi)$;
specifically, $\bpi(C)\leq\exp(-m/576)$.  
Thus the cut $(\bstates_0,\bstates_1)$ is witness to the conductance of the 
single bond flip MC being exponentially small.  
Suppose $d\leq m/12$.  Observe that no $d$-cautious chain can make a transition from  
$\bstates_0 \setminus C$ to $\bstates_1\setminus C$.
Applying Lemma~\ref{lem:torpid},
we therefore obtain.

\begin{theorem}\label{thm:main}
Suppose that $n$, $m$, $G_n$, $\bstates$ and $\bpi$ are as above,
and that $d\leq m/12$. 
Any ergodic Markov chain on state space $\bstates$ 
with stationary distribution~$\bpi$ that is $d$-cautious 
has mixing time $\Omega(\exp(m/576))$.  In particular, the GS Process, which
is $1$-cautious, has mixing time exponential in the number of vertices in $G_n$.
\end{theorem}

It is also natural to consider a ``Swendsen-Wang-style'' Markov chain for sampling from 
$(\ISstates,\omega)$.  Let $I\in\ISstates$ be the current state.  Choose $A$
u.a.r.\ from the set $\{A\in \bstates\mid \chi(I,A)\}$.  Then choose $I'$
u.a.r.\ from the set $\{I'\in\ISstates\mid \chi(I',A)\}$.  The new state is~$I'$.
We can think of this process as a Markov chain on state space $\ISstates\cup\bstates$
with stationary distribution $\frac12\omega$ on~$\ISstates$ and $\frac12\bpi$ 
on~$\bstates$.  (Assume a continuous time process to avoid the obvious periodicity.)
It follows from the earlier analysis that the cut $(\ISstates_0\cup\bstates_0,
\ISstates_1\cup\bstates_1)$ witnesses exponentially small conductance.
To see this, we calculate
the probability in stationarity of observing a transition from 
$\ISstates_0\cup\bstates_0$ to $\ISstates_1\cup\bstates_1$.
There are two 
possibilities:  a transition from $\ISstates_0$ to $\bstates_1$,
or one from $\bstates_0$ to $\ISstates_1$.  The probability of the 
former, we have seen, is 
$\frac12\omega(\ISstates_0)$ times a quantity that is exponentially small in~$n$.
The latter is, by time reversibility, the same as observing, in stationarity,
a transition from $\ISstates_1$ to~$\bstates_0$.  This probability is again exponentially
small in~$n$.  Hence 
the conductance is exponentially small so
the mixing time of the Swendsen-Wang-style Markov chain is exponential in~$n$.

We can also look a little closer, to see what is going on in more 
detail.  Sample a state at random from $(\ISstates,\omega)$,  conditioned on the
event $\omega\in\ISstates_0$, and 
apply a ``half-step'' of the SW-like process to obtain a state $A\in\bstates$.
We know that $A\cap M$ is described by 
a Bernoulli process with success probability~$\frac13$.
Moreover, it is easy to see the remaining edges of~$A$ are Bernoulli with 
success probability~$\frac12$.  Now consider the transition from $A$ to~$I'$.
As in~\cite{GSconf}, view the set $I'\cap U'$ as a $n$-vector~$u'$ over $\mathbb{F}_2$.
Each of the vertices in $V$ that is {\it not\/} incident to an edge of 
$A\cap M$
generates a linear equation, with constant term zero, constraining~$u'$.  
These $\frac23m\approx1.1397 n $ random linear equations 
constrain just $n$ variables;  so with with high probability the only solution is to set
all $n$~variables to~0.  (Equivalently, a random $n\times\frac23m$ 
matrix over $\mathbb{F}_2$ has rank $n$ with high probability.)
In other words, $I'\cap U'=\emptyset$,
except with exponentially small probability, and we find ourselves back
in $\ISstates_0$ again.

\bibliographystyle{plain}
\bibliography{\jobname}
\end{document}